\documentclass[11pt,a4paper]{article}

\usepackage[utf8]{inputenc}
\usepackage{amsfonts}
\usepackage[fleqn]{amsmath}
\usepackage{color}
\usepackage{hyperref}
\usepackage{amsmath,amsthm,amssymb,nccmath}

\usepackage[margin=0.9in]{geometry}
\usepackage{mathtext}
\usepackage[T1,T2A]{fontenc}

\usepackage[matrix,arrow,curve]{xy}
\usepackage{systeme}
\usepackage{ragged2e}
\usepackage{blindtext}
\usepackage{yfonts}
\usepackage{graphicx}  
\usepackage{authblk}

\newtheorem{theorem}{Theorem}

\newtheorem{definition}{Definition}
\newtheorem{remark}{Remark}
\newtheorem{lemma}{Lemma}

\definecolor{dkgreen}{rgb}{0,0.6,0}

\title{Urysohn in action: separating semialgebraic sets by polynomials}

\author[1,2]{Milan Korda}
\author[1]{Jean-Bernard Lasserre}
\author[1]{Alexey Lazarev}
\author[1]{Victor Magron}
\author[3]{Simone Naldi}
 
\affil[1]{LAAS-CNRS, Toulouse, France}
\affil[2]{Faculty of Electrical Engineering, Czech Technical University in Prague}
\affil[3]{Univ. Limoges, CNRS, XLIM, UMR 7252, Limoges, France}
\date{}

\begin{document}

\maketitle



\section{Problem setting}

A classical result from topology called Uryshon's lemma asserts the existence of a continuous separator of two disjoint closed sets in a sufficiently regular topological space. In this work we make a search for this separator constructive and efficient in the context of real algebraic geometry. Namely, given two compact disjoint basic semialgebraic sets $\mathbf{A}=\{x\in \mathbb{R}^n\mid g_i(x)\geq 0\; \forall i = 1 \dots r\}$ and $\mathbf{B}=\{x\in \mathbb{R}^n\mid h_i(x)\geq 0\;\forall i = 1 \dots s \}$ which are contained in an $n$-dimensional box $[-1,1]^n$, we provide an algorithm that computes a separating polynomial $p$ greater than or equal to 1 on $\mathbf{A}$ and less than or equal to 0 on $\mathbf{B}$. This is a challenging problem with many important applications (e.g., classification in machine learning \cite{Kot} or collision avoidance in robotics) and has a long history.   In  \cite{AA99} the authors provide a decision algorithm for the more general separation problem without compactness assumptions.
  In order to obtain a correctness certificate for the separator, another well-renowned approach is to rely on positivity certificates based on sums of squares, such as Putinar certificates (\cite{Putinar93}, \cite[Chapter 2]{Lasserre}) for positive polynomials on basic compact semialgebraic sets. 
  Such certificates have been used to  approximate the volume of a basic semialgebraic set \cite{Henrion,LE} and for the problem of optimal data fitting \cite{LM18}. 
  There, the authors applied hierarchies of semidefinite relaxations coming together with strong convergence guarantees. 
  Related convergence rates have been  obtained \cite{Korda} thanks to the degree bounds for the associated positivity certificates \cite{Nie}. 
 \paragraph{Contributions.}
 Inspired by these latter efforts, we build a Putinar representation of the separating polynomial and estimate the degree of the representation with the help of the recent results \cite{Baldi} which (together with \cite{LS21} and \cite{MM22}) significantly improve the degree estimates in comparison with previous works on that subject \cite{Nie,Schw}.
 We provide a hierarchy of semidefinite programs to compute a separating polynomial for the basic compact semialgebraic sets $\mathbf{A}$ and $\mathbf{B}$ defined above. 
 In addition to that, we estimate the degree of the separating polynomial. 
 Our degree bound is polynomial in the inverse of the euclidean distance between the two sets and singly exponential in the dimension. 
 
 \section{Preliminaries}
 
 In this section we provide the notation and the central results to which we refer in the proof of the main theorem \ref{main result}.

\textbf{Notation.} Let $\mathbb{R}[\mathbf{x}]$ denote a ring of polynomials with $n$ variables, $\mathbb{R}[\mathbf{x}]_{d}$ be a vector space of polynomials with $n$ variables of degree $\leq d$. For $p  \in \mathbb{R}[\mathbf{x}]$ we work with the norm $||p||:=\max_{[-1,1]^n} {p}$, degree of $p$ is $\mathrm{deg}(p)$, $\varepsilon(p):=\frac{\min_\mathbf{A}{p}}{||p||}$. The set of polynomials $\{f_1,\dots, f_t\}$, where $\forall i=1,\dots,t: f_i \in \mathbb{R}[\mathbf{x}]$ is denoted by $\mathbf{f}$ and let $\mathcal{S}(\mathbf{f}):=\{x\in\mathbb{R}^n\mid \forall f\in\mathbf{f}: f(x)\geq 0\}$ be the semialgebraic set defined by $\mathbf{f}$. In these terms the semialgebraic sets in question are $\mathbf{A}=\mathcal{S}(\mathbf{g})$ and $\mathbf{B}=\mathcal{S}(\mathbf{h})$ respectively. A sum of squares of polynomials is denoted by $\Sigma$. Quadratic module generated by $\mathbf{f}$ is the set $\mathcal{Q}(\mathbf{f}):=\{s_{0}+\sum_{i=1}^{t} s_{i} f_{i} \mid s_{i} \in \Sigma\} $. The euclidean distance between two sets $X, Y\subset\mathbb{R}^n$ is denoted by $\mathrm{dist}(X, Y):=\inf\limits_{x\in X,y \in Y}|x-y|$.
 \begin{definition}
 For $l\in \mathbb{N}$ and the set of polynomials $\mathbf{f}$ the truncated quadratic module $\mathcal{Q}_{l}(\mathbf{f})$ of degree $\leq l$ is: 

 $$\mathcal{Q}_{l}(\mathbf{f}):=\left\{s_{0}+\sum_{i=1}^{t} s_{i} f_{i} \mid s_{i} \in \Sigma,\, \operatorname{\mathrm{deg}} (s_{0}) \leq l,\, \operatorname{\mathrm{deg}} (s_{i} f_{i}) \leq l\,   \forall i=1, \ldots, t\right\} .$$
\end{definition} 


\begin{definition}
The modulus of continuity of a function $f:\mathbb{R}^n\rightarrow\mathbb{R}$ is the function $\omega_f(\delta):[0,\infty)\rightarrow \mathbb{R}$, with $\omega_f(\delta):=\sup_{|x-y|<\delta} {|f(x)-f(y)|} $. 
\end{definition}



\begin{theorem}[Multivariate Jackson's Theorem \cite{Newman}]
 \label{Jackson}
 Let $f$ be a continuous function on an $n$-dimensional box $[-1,1]^n$. 
 For each $m\in\mathbb{N}$ there exists a polynomial $p_m$ with $\deg (p_m) = m$ such that for all $x\in [-1,1]^n$ one has
 $$|f(x)-p_m(x)|\leq C \omega_f\left({\frac{n^{3/2}}{m}}\right),$$
 where $C$ is an absolute constant (not depending on the function $f$).
\end{theorem}


\begin{theorem}[\cite{Baldi}]
\label{effective Putinar}
 Assume that $n\geq 2$ and that $f_1, \dots f_t$ satisfy the two  normalization assumptions (1)  $1-||x||_2^2\in\mathcal{Q}(\mathbf{f})$, (2) $\forall i: ||f_i||\leq \frac{1}{2}$. Then every $p$ positive  on $S(\mathbf{f})$ belongs to $\mathcal{Q}_l(\mathbf{f})$ for $l \geq \gamma(n,\mathbf{f}) deg(p)^{3.5nT} \varepsilon(p)^{-2.5nT}$,
 where $\mathfrak{c}$ and $T$ are the Lojasiewicz coefficient and exponent provided in definition 2.4 of \cite{Baldi} and $\gamma(n,\mathbf{f})=n^{3}2^{5nT}r^{n}\mathfrak{c}^{2n}deg(\mathbf{f})^n$.
\end{theorem}
 
 \section{The main result}
 
 We start this section with constructing a continuous separating function $u$ for the sets $\mathbf{A}$ and $\mathbf{B}$. 
 After that we rely on Theorem \ref{Jackson} to provide a uniform polynomial approximation $p$ of $u$. 
 Finally, we show how to compute such a separator $p$ by means of Putinar's representations and obtain an  associated  hierarchy of semidefinite programs (SDP).
 
 \subsection{Explicit continuous separator}
 
 We can build a separating continuous function $u$  explicitly. Consider:
 $$u(x):=2 - 3\frac{\mathrm{dist}(x,\mathbf{A})}{\mathrm{dist}(\mathbf{A},\mathbf{B})}.$$

 \begin{lemma}
 The function $u$ is $L$-lipschitz with $L=\frac{3}{\mathrm{dist}(\mathbf{A},\mathbf{B})}$, $u|_\mathbf{A} = 2$ and $u|_\mathbf{B}\leq -1$.
 \end{lemma}


 \subsection{Uniform approximation of the separator by polynomials }
 Applying Theorem \ref{Jackson} to the function $u$  provides the following bound on the degree $m$ of the approximating polynomial $p_m$:
 $$|u(x)-p_m(x)|\leq C \omega_f\left({\frac{n^{3/2}}{m}}\right),$$ 
 where $C$ is an absolute constant (not depending on  $u$, see \cite{Newman}). For the $L$-lipschitz function $u$ it is clear that: $\omega_u(\delta)\leq L\delta$, and thus we obtain:
 \begin{equation}
 \label{accuracy of uniform appr using Jackson}
     |u(x)-p_m(x)|\leq C L \frac{n^{3/2}}{m},    
 \end{equation}
 where $\deg(p_m)=m$, $L$ is the lipschitz constant of $u$ and $C$ is an absolute constant (not depending on  $u$).

 \subsection{Putinar's representation of the separating polynomial}
 
In this section we provide a Putinar's representation of the separating polynomial $p$ and estimate the degree bound of this representation.  
This separating polynomial is greater than 1 on $\mathbf{A}$ and negative on $\mathbf{B}$.
 
 \begin{remark}
 According to section 2.4.2 of chapter 2 of \cite{Lasserre} it is possible to check membership of a polynomial in $\mathcal{Q}_{l}(\mathbf{g})$ by solving a semidefinite feasibility program.
 \end{remark}
 \begin{theorem}
 \label{main result}
 Assume $\mathfrak{c}$, $T$  and $\gamma(n,\mathbf{f})$ are taken from Theorem \ref{effective Putinar}. Then there exists a separating polynomial of degree $l$ with
 
 
 \begin{equation}
 \label{convergence rate0}
     l \geq \max \{\gamma(n,\mathbf{g}), \gamma(n,\mathbf{h})\} C^{3.5nT}n^{3nT} \left(\frac{6}{\mathrm{dist}(\mathbf{A},\mathbf{B})}\right)^{6nT}.
 \end{equation}
 \end{theorem}
 We can also simplify the result given in formula \ref{convergence rate0}.
 \begin{remark}
If the gradients of the active constraints' defining polynomials from the sets $\mathbf{g}$ and $\mathbf{h}$  are linearly independent at every point of $\mathbf{A} =  \mathcal{S}(\mathbf{g})$ and $\mathbf{B} = \mathcal{S}(\mathbf{h})$, respectively (see precise definition 2.7 in \cite{Baldi}) then $T=1$ and we obtain the bound: 
\begin{equation}
\begin{split}
     l \geq \max \{\gamma(n,\mathbf{g}), \gamma(n,\mathbf{h})\} C^{3.5n}n^{3n} \left(\frac{6}{\mathrm{dist}(\mathbf{A},\mathbf{B})}\right)^{6n}.
\end{split}
 \end{equation}
 \end{remark}
 
 
 

 \section{Practical implementation and discussion}
\begin{figure}[h!]
    \centering
    \includegraphics[scale=0.3]{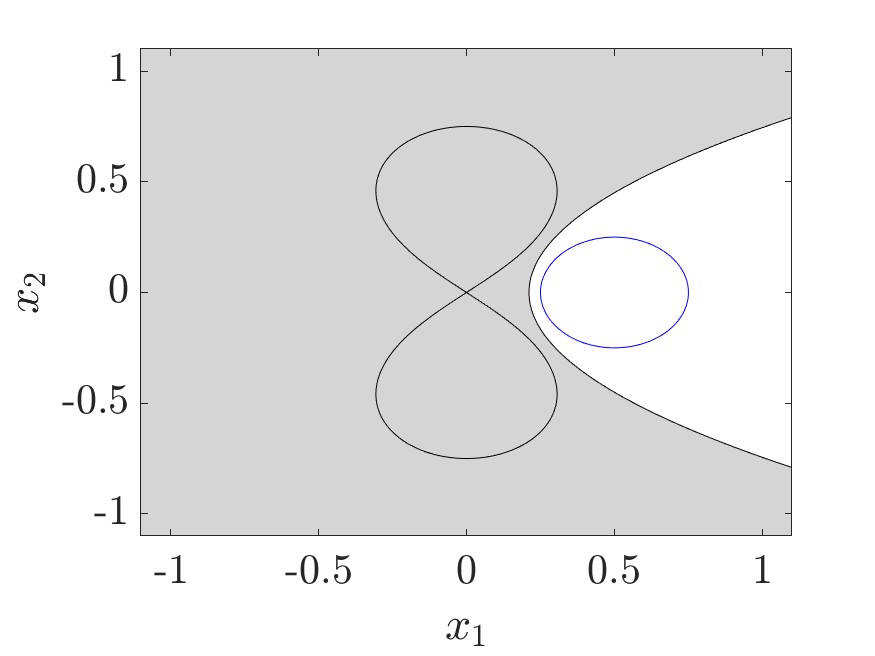}
    \caption{Separating a lemniscate and a circle. The lemniscate and the circle are the boundaries of the compact semialgebraic sets $\mathbf{A}$ and $\mathbf{B}$ respectively. The grey area is the superlevel set of $p$ and the white area is the sublevel set of $p$.
    }\label{fig}
\end{figure}

 Let us consider an example when $\mathbf{A}=\{(x_1,x_2)\in \mathbb{R}^2\mid-\frac{16}{9}(x_1^2+x_2^2)^2+x_2^2-x_1^2\geq 0\}$ and $\mathbf{B}=\{(x_1,x_2)\in \mathbb{R}^2\mid\frac{1}{16}-(x_1^2-\frac{1}{2})^2-x_2^2\geq 0\}$. The levels sets of the separating polynomial $p(x_1,x_2)=1.92876-7.71502x_1+10.96977x_2^2$ with $\mathrm{deg}(p)=2$ are displayed on Figure \ref{fig}. Computations where executed in MATLAB with the usage of Yalmip \cite{Yalmip} and Mosek \cite{Mosek}.\\
 
 \textbf{Perspectives.} Having solved the problem of separating two semialgebraic sets it is immediately possible to construct a polynomial $p$ classifying several semialgebraic sets $\mathbf{A_i}$, i.e.,  $p|_\mathbf{A_i}\in[m_i,M_i]$. 
 Another possible application is to use the described techniques to separate reachable sets of dynamical systems. 
 It could also be interesting to compare the practical efficiency of our current SDP-based framework and  concurrent techniques based on Bernstein polynomials.
 Furthermore, one could try to maximize the separating margin similarly to the way it is done for support-vector machines. 
 Finally, one could exploit sparsity and symmetry patterns arising from the input data in order to improve the scalability of our approach. 



\end{document}